\documentclass{article}
\newcommand{\ml}{l\kern-0.035cm\char39\kern-0.03cm}

\newtheorem{prop}{Proposition 1.}
 
\newtheorem{defi}{Definition 1.}

\newtheorem{ex}{Example}

\newtheorem{prop2}{Proposition 2.}
 
\newtheorem{defi2}{Definition 2.}

\newtheorem{ex2}{Example 2.}

\newtheorem{defi3}{Definition 3.}

\newtheorem{ex3}{Example 3.}
\title {\centerline \bf Map  for simultaneos measurements  for 
a quantum logic}
\author{Olga N\'an\'asiov\'a \footnote
{Supported by grant
VEGA1/7146/20}}
\date{ Slovak Technical University}
\begin{document}
\maketitle
\begin{abstract}
In this paper we will study a function of simultaneus measurements
for quantum events (s-map) which will be  compared with the conditional states 
on an orthomodular lattice as a basic structure for quantum logic. 
We will show the connection between s-map and a conditional
state. Based on the R\'enyi approach to the conditioning, conditional states 
and the induces independence of events with respect to a state are 
discussed. Observe that their relation of independence of events is 
not more symmetric contrary to the standard probabilistic case. 
Some illustrative examples are included.        

\end{abstract}

\vskip 3pc
\centerline{\bf Introduction }
\vskip 2pc

Conditional probability plays a basic role in the classical 
probability theory. Some of the most important areas of the theory such
as martingales, stochastic processes rely heavily of this concept. 
Conditional probabilities on a classical measurable space are studied in
several different ways, but result in equivalent theories.  
The classical probability theory does not decsribe the causality model.

The situation charges when non-standard spaces are considered. 
For example,   it is a well known
 that the set of random events in  quantum mechanics experiments  
is a more general structure
than Boolean algebra. 
 In the quantum logic approach the set of random
events is assumed  to be a orthomodular 
lattice (OML) $L$ .  Such model we 
can find not only in the quantum theory, but for example, in the economics,
biology etc. We will show such such a simple situation in the 
Example 1. 

In this paper we will study  a conditional state  on an
 OML
using  Renyi's approach (or Baye\-sian principle).  This approach   
helps us to define 
independence of events and 
differently from the situation in the classical theory of probability, 
if an event $a$ is independent of an event $b$, then the event $b$
can be dependent on the event $a$ (problem of causality) (\cite{N6},
\cite{N7}).
We will show that we can define a $s$-map (function for simultaneous
measurements on an OML). It can be shown that if we have the 
conditional state we can define the $s$-map and conversely. By using the
$s$-map we can introduce joint distribution also for noncompatible 
observables on an OML. Moreover, if $x$ is an obsevable on $L$
and $ B$ is Boolean sub-algebra of $L$, we can construct an observable 
$z=E(x\vert  B)$, which is a version of conditional expectation of
$x$ but it need not to  be necessarily  compatible with $x$.

\vskip 1pc
\begin{ex}
Assume that there are four objects $(A,U),(A,V),(C,U),(C,V)$
under erxperimental observations, but due to the nature of experimental
 device we are able to identify only one constituent of the pair.
Thus, the possible outcomes of our experiment are $A,C,U,V$, and
if the outcome is (say) $A$ then we do not know whether
it comes from the pair $(A,U)$ or from $(A,V)$. In other words we can
always only one characteristic feature of each object observe:
    $$A=\Pi_1(A,U)=\Pi_1(A,V)\mbox {\hskip 2pc}  C=\Pi_1(C,U)=\Pi_1(C,V)$$
$$U=\Pi_2(A,U)=\Pi_2(C,U)\mbox {\hskip 2pc} V=\Pi_2(A,V)=\Pi_2(C,V)$$
where $\Pi_i$, $i=1,2$  present some  "state" of our system.
In such situation, for example,  we ask about the probability
of $A$ if property $U$ has been detected; equivalently we ask about
the value of $P(A|U)$.
\end{ex}
\vskip 2pc
\centerline{\bf 1. A conditional state on an OML }
\vskip 2pc
In this part we introduce the notions as an OML, a state,
a conditional state and their basic properties.
\begin{defi} Let $L$ be a nonempty set
	endowed  with a partial ordering $\leq $. Let there exists the
	greatest element ($1$) and the smallest element ($0$). Let there
	be defined the operations supremum $(\vee )$, infimum $\wedge $
	(the lattice operations ) and a map $\perp :L\to L$ with the
	following properties:
\begin{itemize}
\item[(i)]    For any $\{a_n\}_{n\in\cal A}\in L$, where
$\cal A\subset N$ are  finite 
	 $$\bigvee_{n\in\cal A} a_n,
                 \bigwedge_{n\in\cal A} a_n\in L.$$

\item[(ii)] For any $a\in L$ $(a^\bot )^\bot =a$.
\item[(iii)] If $a\in L$,  then $a\vee a^\bot =1$.
\item[(iv)]  If $a,b\in L$ such that $a\le b$,  then $b^\bot\le a^\bot $.
\item[(v)]  If $a,b\in L$  such that $a\le b$ then $b=a\vee (a^\bot\wedge b)$ (orthomodular law).
\end{itemize}
	Then
$(L,0,1,\vee ,\wedge ,\perp )$
	is called {\it the orthomodular lattlice} (briefly
{\it  OML}).
\end{defi}
\vskip 1pc
 Let $L$ be OML. Then elements $a,b\in L$  will be called:
\begin{itemize}
\item
{\it orthogonal} ($a\bot b$) iff $a\le b^\bot$;
\item
{\it compatible} ($a\leftrightarrow b$) iff  there exist mutually  
orthogonal elements $a_1,b_1,c\in L$   such that
$$a=a_1\vee c\mbox {\hskip 1pc  and\hskip 1pc } b=b_1\vee c.$$
\end{itemize}
\vskip 1pc
If  $a_i\in L$ for any $i=1,2,3,...$ and
$b\in L$ is such, that $b\leftrightarrow a_i $ for all $i$, then
$ b\leftrightarrow\bigvee_{i=1}^n a_i$ and
$$b\wedge (\bigvee_{i=1}^\infty a_i)=\bigvee_{i=1}^\infty ( a_i\wedge b)$$
(\cite{D2},\cite{P1},\cite{V1}). 
\vskip 1pc A subset $ L_0\subseteq L$ is {\it a
sub-logic of $L$ } if for any $ a\in L_0$ we have $a^\bot\in L_0$
and for any $a,b\in L_0$ $a\vee b\in L_0$.
\vskip 1pc
\begin{defi}
 A map $m:L\to R $  such that
\begin{itemize}
\item[(i)]
    $m(0)=0$ and $m(1)=1$.
\item[(ii)]
    If $a\bot b$ then $m(a\vee b)=m(a)+m(b)$
\end{itemize}
    is called {\it a state} on $L$. If we have
    orthomodular  $\sigma $-lattice and $m$ is
    $\sigma $-additive function, then $m$ will be called a
    {\it $\sigma $-state}.
\end{defi}
\vskip 1pc
\begin{defi}\cite{N7}
Let $L$ be an OML. A subset $L_c\subset L-\{0\} $ is
called {\it a conditional system (CS)  in $L$} {\it ($\sigma$-CS
in $L$ })if the following conditions hold:
\begin{itemize}\item 
If $a,b\in L_c$, then $a\vee b\in L_c$. (If $a_n\in L_c$, for
$n=1,2,...$, then $\bigvee_n a_n\in L_c$.)
\item
If $a,b\in L_c$ and $a< b$,
then $a^\perp\wedge b\in L_c$.
\end{itemize}
\end{defi}
\vskip 1pc
Let $A\subset L$. Then $L_c(A)$ is the smallest CS ($\sigma $-CS
),
which contains the set $A$.

\begin{defi}\cite{N7}
Let $L$ be an OML and $L_c$ be an $\sigma $-CS in $L$.
Let $f:L\times L_c\to [0,1]$. If the function $f$ fulfill the
following conditions: 
\begin{itemize}
\item[(C1)]
for each $a\in L_0$ $f(.,a)$ is a state on $L$;
\item[(C2)]
for each $a\in L_0$ $f(a,a)=1$;
\item[(C3)]
if $\{a_n\}_{n\in\cal A}\in L_0$, where $\cal A\subset N$ and
$a_n$ are mutually orthogonal, then for each $b\in L$
$$f(b,\bigvee_{n\in\cal A}a_n)=\sum_{n\in\cal A}f(a_n,\bigvee_{n\in\cal A}a_n)
f(b,a_n);$$
\end{itemize}
then it is called {\it  conditional state}.
\end{defi}
\vskip 1pc
\begin{prop}\cite{N7}
	Let $L$ be a OML. Let $\{a_i\}_{i=1}^n\in L$, $n\in N$
	where  $a_i\perp a_j$ for $i\neq j$. If
	for any $i$ there exists a state $\alpha_i$,
	such that $\alpha_i(a_i)=1$, then there
exists  $\sigma $-CS such that for any
	${\bf k}=(k_1,k_2,...,k_n)$, where $k_i\in [0;1]$ for
	$i\in\{1,2,...,n\}$ with the property
	$\sum_{i=1}^n k_i=1$,
	there exists a conditional state
	$$f_{\bf k}:L\times L_c\to [0;1],$$
	such that
\begin{enumerate}
\item
	for any $i$ and each $d\in L$
	$f_{\bf k}(d,a_i)=\alpha_i(d);$
\item
	for each $a_i$
	$$f_{\bf k}(a_i,\bigvee_{i=1}^na_i)=k_i;$$
\end{enumerate}
\end{prop}
\vskip 1pc
\begin{defi}\cite{N7}
	 Let $L$ be an OML and $f$ be a  conditional state.
	 Let $b\in L$, $a,c\in L_c$ such that $f(c,a)=1$. Then $b$ is
 independent
	  of $a$ with respect to the state $f(.,c)$ ($b\asymp_{f(.,c)}a$)
	  if $f(b,c)=f(b,a)$.
\end{defi}
\vskip 1pc
 The classical   definition of independency of a probability space
$(\Omega,\cal B,P)$ is a
special case of this definition,
because
\vskip 1pc
\centerline{$P(A|B)=P(A|\Omega )$ if and only if $ P(A\cap B|\Omega )=
P(A|\Omega )P(B|\Omega )$.}
\vskip 2pc
 
 If $L_c$ be CS and  $f:L\times L_c\to [0,1]$
 is a conditional state, then (\cite{N7})
 \begin{itemize}\item[(i)]
 Let $a^\perp,a,c\in L_c$, $b\in L$ and $f(c,a)=f(c,a^\perp )=1$. Then
 $b\asymp_{f(.,c)} a$  if and only if $b\asymp_{f(.,c)} a^\perp$.
 \item[(ii)]
  Let $a,c\in L_c$, $b\in L$ and $f(c,a)=1$. Then
 $b\asymp_{f(.,c)} a$  if and only if $b^\perp \asymp_{f(.,c)} a$.
 \item[(iii)]
 Let $a,c,b\in L_c$, $b\leftrightarrow a$ and $f(c,a)=f(c,b)=1$. Then
 $b\asymp_{f(.,c)} a$  if and only if $a\asymp_{f(.,c)} b$.
 \end{itemize}

\vskip 5pc
\centerline {\bf 2. Function for simultaneous measurement ({\it s-map})}
\vskip 1pc
\begin{defi2}
Let $L$ be an OML. The map $p:L\times L\to [0,1]$
will be called {\it s-map} if the following conditions  hold:
\begin{itemize} 
\item[(s1)]
$p(1,1)=1$;
\item[(s2)]
if $a\perp b$, then $p(a,b)=0$;
\item[(s3)]
if  $a\perp b$, then for any $c\in L$,
$$p(a\vee b,c)=p(a,c)+p(b,c)$$
$$p(c,a\vee b)=p(c,a)+p(c,b)$$.
\end{itemize}
\end{defi2}
\vskip 1pc
\begin{prop2}
Let $L$ be an OML and let $p$ be a s-map.
Let $a,b,c\in L$, then 
\begin{itemize}
\item[1.]
if $a\leftrightarrow b$ , then $p(a,b)=p(a\wedge b,a\wedge b)=p(b,a)$;
\item[2.] 
if  $a\le b$, then  $p(a,b)=p(a,a)$;
\item[3.]
if $a\le b$, then $p(a,c)\le p(b,c)$;
\item[4.]
 $p(a,b)\le p(b,b)$;
\item[5.]
if  $\nu (b)=p(b,b)$, then $\nu $ is a state on $L$.
\end{itemize}
\end{prop2}
Proof.
{\it (1)} If $a\leftrightarrow b$, then
$a=(a\wedge b)\vee (a\wedge b^\perp)$ and
$b=(b\wedge a)\vee (b\wedge a^\perp)$. Hence
$$p(a,b)=p((a\wedge b)\vee (a\wedge b^\perp),b)=$$
$$=p(a\wedge b,b)+p(a\wedge b^\perp ,b)=p(a\wedge b, b).$$
Analogously
$$p(a\wedge b, b)=p(a\wedge b,(b\wedge a)\vee (b\wedge a^\perp))=$$
$$=p(b\wedge a,b\wedge a)+p(b\wedge a,b\wedge a^\perp)=p(b\wedge a,b\wedge a)
. $$
Hence
$$p(a,b)=p(a\wedge b,a\wedge b).$$
{\it (2)} If $a\le b$, then $a\leftrightarrow b$. Hence
 $$p(a,b)=p(a,a\wedge b)=p(a,a).$$
{\it (3)} If $a\le b$, then $b=a\vee (a^\perp\wedge b)$. Hence
$$\begin{array}{clcr} p(b,c)&=p(a\vee (a^\perp\wedge b),c)\\
               &=p(a,c)+p(a^\perp\wedge b,a)\ge p(a,c)
\end{array}
$$
{\it (4)} From (3) and {\it (2)} it follows 
$$p(b,b)=p(1,b)\ge p(a,b).$$ 
 Hence we get 
$$p(b,b)\ge p(a,b)\mbox{ \hskip 2pc for each\hskip 2pc } a,b\in L.$$
\newline
{\it (5)} Let $\nu :L \to [0,1]$, such that
$\nu (b)=p(b,b)$. Then
$$\nu (0)=p(0,0)=0.$$
Let $a\perp b$, then
$$\nu (a\vee b)=p(a\vee b,a\vee b)=p(a,a\vee b)+p(b, a\vee b)=$$
$$=p(a,a)+p(a,b)+p(b,a)+p(b,b)=p(a,a)+p(b,b)=\nu (a)+\nu (b).$$
From the definition we have that $\nu (1)=p(1,1)=1$.
From this it follows that $\nu $
is a state on $L$.
\newline
(Q.E.D.)
\vskip 1pc
\begin{prop2}
Let $L$ be an OML, let there be a s-map $p$.
Then there exists  a  conditional state $f_p$, such that
$$p(a,b)=f_p(a,b)f_p(b,1).$$

Let $L$ be an OML and let $L_c=L-\{0\}$. If
$f:L\times L_c\to [0,1]$
is a conditional state, then  there exists a s-map $p_f:L\times L\to [0,1]$.
\end{prop2}
Proof.
Let $p$ be a s-map. Let $L_c=\{b\in L;\quad p(b,b)\neq 0\}$.
Let $f_p:L\times L_c\to  R$ such that
$$f_p(.,b)=\frac{p(.,b)}{p(b,b)}. $$
From the Proposition 2.1 (3) it follows that
for any $a\in L$ and $b\in L_c$ $f_p(a,b)\in [0,1]$.
Moreover
$$f_p(0,b)=0\mbox{\hskip 2pc and\hskip 2pc } f_p(1,b)=\frac{p(1,b)}{p(b,b)}=
\frac{p(b,b)}{p(b,b)}=1$$
and also $f_p(b,b)=1$.
Let $c,a\in L$ and let $a\perp c$. Then 
$$ 
f_p(a\vee c,b)=\frac{p(a\vee c,b)}{p(b,b)}=\frac{p(a,b)+p(c,b)}{p(b,b)}
=f_p(a,b)+f_p(c,b).$$ 
It means that for any $b\in L_c$ is $f_p(.,b)$ a state on $L$.

Let $b_i\in L_c$, $i=1,2,...,n$ be mutually orthogonal elements.
Then for any $a\in L$
$$
 f_p(a,\bigvee_{i=1}^nb_i)=\frac{p(a,\vee_ib_i)}
{p(\vee_ib_i,\vee_ib_i)}
=\sum_{i=1}^n\frac{p(a,b_i)}{p(\vee_ib_i,\vee_ib_i)}
=\sum_{i=1}^n\frac{p(b_i,\vee_ib_i)}{p(\vee_ib_i,\vee_ib_i)}
\frac{p(a,b_i)}{p(b_i,\vee_ib_i)}$$
$$= \sum_{i=1}^n\frac{p(b_i,\vee_ib_i)}{p(\vee_ib_i,\vee_ib_i)}
\frac{p(a,b_i)}{p(b_i,b_i)}=\sum_{i=1}^nf_p(b_i,\vee_ib_i)f(a,b_i).
$$

From this it follows that $f_p$ is the conditional state.

Now we can compute
$$ f_p(a,b)f_p(b,1)=\frac{p(a,b)}{p(b,b)}\frac{p(b,1)}{p(1,1)}.$$
From the properties of s-map we have $p(b,1)=p(b,b)$ and $p(1,1)=1$.
Hence $f_p(a,b)f_p(b,1)=p(a,b)$.

Let $f$ be a conditional state and let $L_0=\{b\in L_c;f(b,1)\ne 0\}$.
Let 
$$p_f: L\times L\to [0,1]$$ be defined in the following way:

$$p_f(a,b)= \left\{\begin{array}{l}f(a,b)f(b,1),\quad b\in L_0\\
0,\quad\quad\quad\quad\quad\quad  b\notin L_0\end{array}\right.$$

(s1) Because $1\in L_0$ and $f$ is a conditional state, then
$$p_f(1,1)=f(1,1)f(1,1)=1.$$

(s2) Let $a,b\in L$ and $a\perp b$. If $b\in L_0$, then
$p_f(a,b)=f(a,b)f(b,1)$.
Because $a\le b^\perp$, then $f(a,b)=0$. Hence $p_f(a,b)=0$.
If $b\notin L_0$, then then $p_f(a,b)=0$. Hence for $a\perp b$
$p_f(a,b)=0$.

(s3)
Let $a,b,c\in L$, $a\perp b$. We have to show that 
\begin{equation} p_f(a\vee b,c)=p_f(a,c)+p_f(b,c)\end{equation}
and
\begin{equation} p_f(c,a\vee b)=p_f(c,a)+p_f(c,b).\end{equation}

(1) If $c\in L_0$, then

$$\begin{array}{clcr}p_f(a\vee b,c)&=f(a\vee b,c)f(c,1)\\
&=f(a,c)f(c,1)+f(b,c)f(c,1)\\
&=p_f(a,c)+p_f(a,c).
\end{array} $$
If $c\notin L_0$, then $p_f(a\vee b,c)=p_f(a,c)=p_f(b,c)=0$.
Hence $$p_f(a\vee b,c)=p_f(a,c)+p_f(b,c).$$

(2) In this case we have to verify for (b) the following three situations:
\begin{itemize} \item[(i) ] $a,b\in L_0$;
\item[(ii)]  $a\in L_0$, $b\notin L_0$;
\item[(iii)] $a,b\notin L_0$.
\end{itemize}

(i) If $a,b\in L_0$, then
$$\begin{array}{clcr} p_f(c,a\vee b)&=f(c,a\vee b)f(a\vee b,1)\\
&=(f(a,a\vee b)f(c,a)+f(b,a\vee b)f(c,b))f(a\vee b,1)\\
&=f(c,a)f(a,a\vee b)f(a\vee b,1)+f(c,b)f(b,a\vee b)f(a\vee b,1).
\end{array}
$$
From the definition  of the function $f$ we get
$$\begin{array}{clcr} f(a,1)&=f(a,a\vee b)f(a\vee b,1)+
f(a,(a\vee b)^\perp)f((a\vee b)^\perp,1)\\
&=f(a,a\vee b)f(a\vee b,1)+0.\end{array} $$
Also
$$f(b,a\vee b)f(a\vee b,1)=f(b,1).$$
Then
$$\begin{array}{clcr} p_f(c,a\vee b)&=
f(c,a)f(a,a\vee b)f(a\vee b,1)+f(c,b)f(b,a\vee b)f(a\vee b,1)\\
&=f(c,a)f(a,1)+f(c,b)(f(b,1)\\
&=p_f(c,a)+p_f(c,b).\end{array}$$

 (ii) If $a\in L_0$ and $b\notin L_0$ and $a\vee b\in L_0$,  then
 from the definition a map $p_f$ it follows $p_f(c,b)=0$. From this it
 follows that it is enought to show
 $$p_f(c,a\vee b)=p_f(c,a).$$
 But
 $$p_f(c,a\vee b)=f(c,a\vee b)f(a\vee b,1)$$
 and
 $$p_f(c,a)=f(c,a)f(a,1).$$
  
 Because $f(b,1)=0$, then
 $$f(a\vee b,1)=f(a,1)+f(b,1)=f(a,1).$$ 
 On the other hand
 $$0=f(b,1)=f(a\vee  b,1)f(b,a\vee b)+
 f((a\vee b)^\perp ,1)f(b,(a\vee b)^\perp ).$$
 Because $f(b,(a\vee b)^\perp )=0$, then we have
 $$0=f(a\vee b,1)f(b,a\vee b).$$
 But  $f(a\vee b,1)\ne 0$ and hence
 $$f(b,a\vee b)=0$$
 and so
 $$1=f(a\vee b,a\vee b)=f(a,a\vee b)+f(b,a\vee b)=f(a,a\vee b).$$
Therefore
$$f(c,a\vee b)=f(a,a\vee b)f(c,a)+f(b,a\vee b)f(c,b)=f(c,a).$$
Hence
$$p_f(c,a\vee b)=f(c,a\vee b)f(a\vee b,1)=f(c,a)f(a,1)=p_f(c,a).$$
(iii) If $a,b\notin L_0$, then  $f(a,1)=f(b,1)=0$. From this it follows that
$f(a\vee b,1)=0$ and so $a\vee b\notin L_0$.
Hence for any $c\in L$
$$0=p_f(c,a\vee b)=p_f(c,a)+p_f(c,b).$$
 
Therefore  $p_f$ is s-map.
(Q.E.D.)
  \vskip 1pc
 \begin{prop2}
 Let $L$ be an OML.
 \begin{itemize}\item[(a)]
 If $f$ is a conditional state, then
 $b\asymp_{f(.,1)} a$  iff $p_f(b,a)=p_f(a,a)p_f(b,b)$, where
 $p_f$ is the s-map generated by $f$.
 \item[(b)]
 Let $p$ be a s-map. Then $b\asymp_{f_p(.,1)} a$
 iff $p(b,a)=p(a,a)p(b,b)$, where
 $f_p$ is the conditional state generated by the s-map $p$.
 \end{itemize}
 \end{prop2}
 Proof.

 (a) Let $b\asymp_{f(.,1)}a$. It means that 
 $f(b,a)=f(b,1)$.
 Let $f(b,1)\ne 0$ and $f(a,1)\ne 0$. 
 From the previous proposition we have that
 $$p_f(b,a)=f(b,a)f(a,1)=f(b,1)f(a,1).$$
 But
 $$p_f(d,d)=f(d,d)f(d,1)=f(d,1)$$
 and  hence
 $$p_f(b,a)=p_f(b,b)p_f(a,a).$$

 Let $f(b,1)=0$ and $f(a,1)\ne 0$. From this it follows 
 that $p_f(b,b)=f(b,1)=0$. On the other hand
 $$f(b,1)=f(a,1)f(b,a)+f(a^\perp ,1)f(b,a^\perp )=0.$$
 Therefore
 $f(b,a)=0$ and hence $p_f(b,a)=0$. It means that in this case 
 $p_f(b,a)=p_f(b,b)p_f(a,a)$.
 
 Let $f(b,1)=f(a,1)=0$. From this it follows that $f(a,1)=p_f(a,a)=0=p_f(b,b)$
 and so $p_f(a,a)p_f(b,b)=0$. On the other hand $p_f(b,a)=f(b,a)f(a,1)=0$.
 It means
 
\begin{equation} b\asymp_{f(.,1)}a\mbox {\hskip 1pc implies\hskip 1pc } p_f(b,a)=p_f(a,a)p_f(b,b)
 .\end{equation}
 
 If $p_f(b,a)=p_f(a,a)p_f(b,b)$, then $p_f(b,a)=f(a,1)f(b,1)$.
 It means that
 $$p_f(b,a)=f(b,a)f(a,1)=f(b,1)f(a,1).$$
 From this  it follows
 $$f(b,1)=f(b,a),$$
 and so
 $$b\asymp_{f(.,1)}a.$$

 (b) Let $p$ be a s-map and $L_c=\{d\in L;\quad p(d,d)\ne 0\}$. Let
 $f_p:L\times L_c\to [0;1]$ be the conditional state
 defined
$$f_p(b,a)=\frac{p(b,a)}{p(a,a)}.$$ Let
 $b\asymp_{f_p(.,1)}a$. It means that $f_p(b,a)=f_p(b,1)$. Hence
 $$f_p(b,a)=\frac{p(b,a)}{p(a,a)}$$
 and
 $$f_p(b,1)=\frac{p(b,1)}{p(1,1)}=p(b,b).$$
 Hence
 $$\frac{p(b,a)}{p(a,a)}=p(b,b)$$
 and so
 $$p(b,a)=p(a,a)p(b,b).$$

 On the other hand, if $p(a,b)=p(a,a)p(b,b)$, then
 $$\begin{array}{clcr} f_p(b,a)&=\frac{p(b,a)}{p(a,a)}=\frac{p(a,a)p(b,b)}{p(a,a)}\\
 &=p(b,b)=p(b,1)\\&=\frac{p(b,1)}{p(1,1)}
  =f_p(b,1).\end{array}$$
 It means $b\asymp_{f_p(.,1)}a$.

 (Q.E.D.)

\vskip 1pc
\begin{ex2}
Let  $L=\{a,a^\perp ,b,b^\perp , 0,1\}$. It is clear that
$L$ is an OML. Let $f(s,t)$ is defined by the following way:

$$
\vbox{\offinterlineskip
\halign{
\strut\vrule $#$ & \vrule $#$ &\vrule $#$
&\vrule $#$ &\vrule $#$&\vrule $#$\hfil\vrule \cr
\noalign{\hrule}
\quad s/t   \quad &\quad a\quad &\quad a^\perp\quad &\quad b   \quad &\quad b^\perp \quad&\quad 1  \quad\cr
\noalign{\hrule}
\quad a     \quad &\quad 1\quad &\quad 0 \quad &\quad 0.4 \quad &\quad 0.4    \quad &\quad 0.4\quad\cr
\noalign{\hrule}
\quad a^\perp\quad &\quad 0\quad &\quad 1      \quad &\quad 0.6\quad & \quad 0.6   \quad &\quad 0.6\quad\cr
\noalign{\hrule}
\quad b\quad &\quad 0.2\quad &\quad11/30\quad &\quad 1\quad &\quad 0\quad &\quad 0.3\quad\cr
\noalign{\hrule}
\quad b^\perp\quad &\quad 0.8\quad &\quad 19/30\quad &\quad 0\quad  & \quad 1\quad&\quad 0.7\quad \cr
\noalign{\hrule}
}}$$
From $f$ this we can compute $p_f(s,t)$ . Then we get:
$$
\vbox{\offinterlineskip
\halign{
\strut\vrule $#$ & \vrule $#$ &\vrule $#$
&\vrule $#$ &\vrule $#$\hfil\vrule \cr
\noalign{\hrule}
\quad s/t \quad &\quad a\quad &\quad a^\perp\quad &\quad b   \quad &\quad b^\perp \quad\cr
\noalign{\hrule}
\quad a\quad &\quad 0.4 \quad &\quad 0\quad &\quad 0.12 \quad &\quad 0.28   \quad \cr
\noalign{\hrule}
\quad a^\perp\quad &\quad 0\quad &\quad 0.6\quad &\quad 0.18 \quad & \quad 0.42 \quad \cr
\noalign{\hrule}
\quad b\quad &\quad 0.08\quad &\quad 0.22\quad &\quad 0.3\quad &\quad 0\quad \cr
\noalign{\hrule}
\quad b^\perp\quad &\quad 0.32\quad &\quad 0.38\quad &\quad 0\quad  & \quad 0.7\quad \cr
\noalign{\hrule}
}}$$
We can see that $p_f(a,b)=p_f(a,a)p_f(b,b)$, but $p_f(b,a)\ne p_f(b,b)p_f(a,a)$.
\end{ex2}
\vskip 2pc
 \centerline{\bf  3. On observables }
  \vskip 1pc
  Let $\cal B(R)$ be $\sigma$-algebra of Borel sets.
A  $\sigma $-homomorphism   $x:\cal B(R)\to L$ is called an observable on $L$.
   If $x$ is an observable, then $R(x):=\{x(E);\quad E\in\cal F \}$
   is called  range of the observable $x$. It is clear that $R(x)$ is
Boolean $\sigma $-algebra [Var].  Let us denote $\nu (b)=p(b,b)$ for $b\in L$.
   \begin{defi3}
   Let $L$ be a $\sigma $-OML and $p:L\times L\to [0;1]$
   be a s-map.
   Let $x,y$ be some observables on $L$.  Then a map
   $p_{x,y}:\cal B(R)\times \cal B(R)\to [0,1]$, such that
   $$p_{x,y}(E,F)=p(x(E),y(F)),$$
   is called a joint distribution for the observables $x$ and $y$.
   \end{defi3}
   \vskip 1pc

    If $F_{x,y}(r,s)=p(x(-\infty ,r),y(-\infty, s))$, then 
the  function $F_{x,y}$ is {\it the distribution  function of 
the observables} 
$x,y$. It is clear that   for $r_1\le r_1$, then 
 $F_{x,y}(r_1 ,s)\le F_{x,y}(r_2,s))$.
  
	If $x$ is an observable on $L$ and 
$m$ is a state on $L$, then $m_x(E)$, $E\in\cal B(R)$
 is probability distribution for $x$ and  
$$m(x)=\int_R\lambda m_x(d\lambda )$$ 
 is called the expectation of $x$ in the state $m$, 
 if the integral on the right 
side exists. 
\vskip 1pc
\begin{defi3} Let $x$ be an obsevable on $L$ and $ B$ be 
a Boolean sub-algebra of $L$ and $f$ be conditional state on $L$
such that $L_c=L-\{0\}$.
Then the observable $z$  will be called 
{\it a conditional 
expectation of} $x$ {\it with respect to } $ B$ {\it in the 
state } $f(.,1)$ iff for any $b\in B-\{0\}$ 
$$f(x,b)=f(z,b).$$ 
We will denote $z:=E_f(x\vert B)$. 
\end{defi3}
\vskip 1pc
It is clear that if 
$L$ be a Boolean algebra, then $E_f(x\vert B)$ is 
known the conditional expectation. The expectation of $x$ in 
the state $m$ have been studied in many papers \cite{G1}-\cite{G8},\cite{P1},
\cite{D2},\cite{N1},\cite{N3}-\cite{N5}, etc. 
In the 
end we show that such the conditional expectation can exist 
on $L$. 
\vskip 1pc
\begin{ex3}
Let  $L$ be the same as in the Example 2.1. 
We have the set 
$$\{f(.,a),f(.,a^\bot ), f(.,b) ,f(.,b^\bot ),f(.,1)\}$$ 
of  states  and 
 $B_d=\{0,1,d,d^\bot\}$, where $d\in L$. Let $x,z $ be  observales
on $L$ such that $R(x)=B_a$, and $R(z)= B_b$. It is easy to 
see, that $x$ is not compatible with $z$. Let
$$x(r_1)=a\quad x(r_2)=a^\bot$$
$$z(s_1)=b\quad z(s_2)=b^\bot $$
for  $r_1,r_2,s_1,s_2\in R$.

If $z=E_f(x\vert  B)$, then  
$$f(x,b)=f(z,b),\quad f(x,b^\bot)=f(z,b^\bot),\quad f(x,1)=f(z,1).$$
From the definition of the  expectation of an observable  we have
$$\begin{array}{clcr}
f(x,1)&=r_1f(a,1)+r_2f(a^\bot ,1)=f(z,1)
\\&=s_1f(b,1)+s_2f(b^\bot ,1 ),\\
f(x,b)&=r_1f(a,b)+r_2f(a^\bot ,b)
=f(z,b)\\&=s_1f(b,b)+s_2f(b^\bot ,b )=s_1,\\
f(x,b^\bot )&=r_1f(a,b^\bot )+r_2f(a^\bot ,b^\bot )=f(z,b^\bot )
\\&=
s_1f(b,b^\bot )+s_2f(b^\bot ,b^\bot )=s_2.
\end{array}$$
 Let $s_1\ne s_2$. If we put 
$$s_1=r_1f(a,b)+r_2f(a^\bot ,b)$$
and
$$s_2=r_1f(a,b^\bot )+r_2f(a^\bot ,b^\bot ),$$
then
$$\begin{array}{clcr} f(z,1)&=s_1f(b,1)+s_2f(b^\bot,1)\\
&=[r_1f(a,b)+r_2f(a^\bot ,b)]f(b,1)+[r_1f(a,b^\bot )
+r_2f(a^\bot ,b^\bot )]f(b^\bot ,1)\\
&=r_1[f(a,b)f(b,1)+f(a,b^\bot )f(b^\bot ,1)]
\\ &\quad +r_2[f(a^\bot ,b)f(b,1)+f(a^\bot ,b^\bot )f(b^\bot ,1)]\\
&=r_1f(a,1)+r_2f(a^\bot ,1)
=f(x,1).\end{array} $$
From this it follows that $z=E_f(x\vert B)$.

If $a\asymp_{f(.,1)} b$, then $f(a,b)=f(a,1)=f(a,b^\bot)$.
From the definition of the  expectation of an observable  we have
$$\begin{array}{clcr}
f(x,b)&=r_1f(a,1)+r_2f(a^\bot ,1)
=f(z,b)=f(z,1)=s,\\
f(x,b^\bot )&=r_1f(a,1)+r_2f(a^\bot ,1 )=f(z,b^\bot )=f(z,1)
=s\\
f(x,1)&=r_1f(a,1)+r_2f(a^\bot ,1)=f(z,1)
\\&=s_1f(b,1)+s_2f(b^\bot ,1 )=s(f(b,1)+f(b^\bot,1))=s.
\end{array}$$
Therefore 
$$f(x,1)=f(x,b)=f(x,b^\bot)=f(z,1)=s,$$
then $R(z)=\{0,1\}\subset B_b$, $z(s)=1$ and  moreover $z=E_f(x\vert B_b)$.

The joint distribution for the observables $x,y$ is given in the table 2.
The second and the third collumms are  $p_{x,y}$ and the fourth and the 
fifth collumms are $p_{y,x}$.

If $R(x)=B_a$ and $x(1)=a$, $x(2)=a^\bot$, then 
$$f(x,1)= f(x,b)=f(x,b^\bot)=1.6.$$
Let $z:=E_f(x\vert B_b)$. Hence
$$f(x,1)=f(z,1)=f(z,b)=f(z,b^\bot)=1.6.$$
Therefore $E_f(x\vert B_b)(1.6)=1$. (In the Example 2.1 
for any $d\in B_b-\{0\}$ and any $c\in B_a$ $c\asymp_{f(.,1)}d$.)

On the other hand, let    $R(y)=B_b$, $y(1)=b$, $y(2)=b^\bot$ and
$w:=E_f(y\vert B_a)$.
Hence  
$$f(y,1)=1.7=0.4w_1+0.6w_2$$
$$f(y,a)=1.8=w_1,\quad f(y,a^\bot)=\frac{49}{30}=w_2$$
and so 
$$E_f(y\vert B_a)(1.8)=a,\quad E_f(y\vert B_a)(\frac{49}{30})=a^\bot .$$

\end{ex3}

\noindent {\bf Author:}  O\v lga N\'an\'asiov\'a,\\
Department of Mathematics and Descriptive Geometry,\\ Faculty of Civil
Engineering,\\ Slovak University of Technology,\\
Radlinsk\'eho 11, \\813 68 Bratislava, \\Slovakia

\noindent e-mail: olga@vox.svf.stuba.sk

\end{document}